\documentclass{article}

\usepackage{amsmath}
\usepackage{amssymb}  
\usepackage{epsfig}
\usepackage{amsthm}   
\usepackage[mathcal]{eucal}
\usepackage{url} 
\newcommand{\F}{\mathcal{F}}

\newcommand{\R}{\mathbb{R}}
\newcommand{\BR}{\bar{\mathbb{R}}}

\DeclareMathOperator{\clconv}{cl\,conv}

\newcommand{\inner}[2]{\langle{#1},{#2}\rangle}

\newcommand{\norm}[1]{\|#1\|}
\newcommand{\normq}[1]{ {\|#1\|}^2 }

\newcommand{\tos}{\rightrightarrows} 

\newtheorem{theorem}{Theorem}[section]

\title{A new old class of maximal monotone operators}

\author{M. Marques Alves\thanks{IMPA, Estrada Dona Castorina 110, 22460-320
    Rio de Janeiro, Brazil
   ({\tt maicon@impa.br})}\hspace{.5em}\thanks{Partially supported by Brazilian CNPq
    scholarship 140525/2005-0.}
  \and
    B. F. Svaiter\thanks{ IMPA, Estrada Dona Castorina 110, 22460-320 Rio de
    Janeiro, Brazil ({\tt benar@impa.br}) }\hspace{.5em}
    \thanks{Partially supported by CNPq
    grants 300755/2005-8, 475647/2006-8 and by PRONEX-Optimization}
}

\date{}

\begin{document}

\maketitle

\begin{abstract}
  In a recent paper in Journal of Convex Analysis the authors studied,
  in non-reflexive Banach spaces, a class of maximal monotone
  operators, characterized by the existence of a function in
  Fitzpatrick's family of the operator which conjugate is above the
  duality product.
  This property was used to prove that such operators satisfies a
  restricted version of Br\o ndsted-Rockafellar property.

  In this work we will prove that if a single Fitzpatrick function of
  a maximal monotone operator has a conjugate above the duality
  product, then all Fitzpatrick function of the operator have a
  conjugate above the duality product.
  As a consequence, the family of maximal monotone operators with this
  property is just the class NI, previously defined and studied by
  Simons.

  We will also prove that an auxiliary condition used by the authors to prove
  the restricted Br\o ndsted-Rockafellar property is equivalent to the
  assumption of the conjugate of the Fitzpatrick function to majorize
  the duality product.
%
  \\
\end{abstract}

\pagestyle{plain}


\section{Introduction}

Let $X$ be a real Banach space.  We use the notation $X^*$ for the
topological dual of $X$ and 
$\inner{\cdot}{\cdot}$ for the duality product
in $X\times X^*$:
\[
\inner{x}{x^*}=x^*(x).
\]
Whenever necessary, we will identify $X$ with its image under the
canonical injection of $X$ into $X^{**}$. 
To simplify the notation, from now on $\pi$ and
$\pi_*$ stands for the duality product in $X\times X^*$ and $X^*\times
X^{**}$ respectively:
\begin{align}
 \nonumber
  &\pi:X\times X^*\to\R, &&\pi_*:X^*\times X^{**}\to\R\\
  \label{eq:df.pps}
  &\pi(x,x^*)=\inner{x}{x^*},&  &\pi_*(x^*,x^{**})=\inner{x^*}{x^{**}}.
\end{align}
The \emph{indicator function} of $A\subset X$ is $\delta_A:X\to\BR$,
\[\delta_A(x):=
\begin{cases}
  0,& x\in A\\
  \infty,& \mbox{ otherwise.}
\end{cases}
\]
For $f:X\to\BR$, the lower semicontinuous convex closure of $f$ is
$\clconv f:X\to\BR$, the largest lower semicontinuous convex function
majorized by $f$. The conjugate of $f$ is $f^*:X^*\to \BR$,
\[
 f^*(x^*)=\sup_{x\in X} \inner{x}{x^*}-f(x).
\]
It is trivial to check that
$ f^*=(\clconv f)^*$.

 A point point-to-set
operator $T:X\tos X^*$ is a relation on $X\times X^*$:
\[ 
 T\subset X\times X^* 
\]
and $x^*\in T(x)$ means $(x,x^*)\in T$.
An operator $T:X\tos X^*$ is {\it monotone} if
\[
\inner{x-y}{x^*-y^*}\geq 0,\forall (x,x^*),(y,y^*)\in T.
\]
and it is {\it maximal monotone}  if it
is monotone and maximal (with respect to the inclusion) in the family
of monotone operators of $X$ in $X^*$.
Maximal monotone operators in Banach spaces arises, for example, in
the study of PDE's, equilibrium problems and calculus of variations.

Given a maximal monotone operator $T:X\tos X^*$, Fitzpatrick
defined~\cite{Fitz88} the family $\mathcal{F}_T$ as those convex, lower
semicontinuous functions in $X\times X^*$ which are bounded bellow by
the duality product and coincides with it at $T$:
\begin{equation}
  \label{eq:def.ft}
  \F_T=\left\{ h\in \BR^{X\times X^*}
    \left|
      \begin{array}{ll}
        h\mbox{ is convex and lower semicontinuous}\\
        \inner{x}{x^*}\leq h(x,x^*),\quad \forall (x,x^*)\in X\times X^*\\
        (x,x^*)\in T 
        \Rightarrow 
        h(x,x^*) = \inner{x}{x^*}
      \end{array}
    \right.
  \right\}.
\end{equation}
Fitzpatrick found an explicit formula for the minimal element of
$\mathcal{F}_T$, from  now on \emph{Fitzpatrick function} of $T$,
$\varphi_T:X\times X^*\to\BR$
\begin{equation}
  \label{eq:def.f.fitz}
  \varphi_T(x,x^*)=\sup_{(y,y^*)\in T} \inner{x}{y^*}+\inner{y}{x^*}-
    \inner{y^*}{y}.
\end{equation}
Moreover, he also proved that if $h\in \F_T$ then $h$ represents $T$
in the following sense:
\[
 (x,x^*)\in T\iff h(x,x^*)=\inner{x}{x^*}.
\]
Note that
\[ \varphi_T(x,x^*)=(\pi+\delta_T)^*(x^*,x).\]
The supremum of Fitzpatrick family is the $\mathcal{S}$-function,
defined and studied by Burachik and Svaiter in~\cite{BuSvSet02},
$\mathcal{S}_T:X\times X^*\to\BR$
\[
\mathcal{S}_T(x,x^*)=\sup  \left\{ h(x,x^*)\;\left|\;
    \begin{array}{l}
  h: X\times X^*\to \BR  \mbox{  convex lower semicontinuous}\\
  h(x,x^*)\leq \inner{x}{x^*}, \quad\forall (x,x^*)\in T
    \end{array}
  \right\}\right.
\]
or, equivalently (see~\cite[Eq.(35)]{BuSvSet02}, \cite[Eq. 29]{BuSvIMPA01})
\begin{equation}
  \label{eq:def.sf}
  \mathcal{S}_T=\clconv (\pi+\delta_T).
\end{equation}
Some authors~\cite{BorJCA06,VosSet06,BorProc07} attribute the
$\mathcal{S}$-function to~\cite{PenRelv04} although \cite{PenRelv04}
was submitted after the publication of~\cite{BuSvSet02}.
Moreover, the content of~\cite{BuSvSet02}, and specifically the $\mathcal{S}$
function, was presented on Erice workshop
on July 2001, by R. S. Burachik~\cite{BuErice01b}.  A list of the
talks of this congress, which includes~\cite{PenErice01}, is available
on the www\footnote{
  \url{http://www.polyu.edu.hk/~ama/events/conference/EriceItaly-OCA2001/Abstract.html}}.
It shall also be noted that~\cite{BuSvIMPA01}, the preprint
of~\cite{BuSvSet02}, was published ( and available on www) at IMPA
preprint server in August 2001.

Burachik and Svaiter defined~\cite{BuSvSet02}, for $h:X\times X^*\to \BR$,
\[
 \mathcal{J}h:X\times X^*\to\BR,
 \qquad  \mathcal{J}h(x,x^*)=h^*(x^*,x)
\]
and proved that if  $T$ is maximal monotone, then
 $\mathcal{J}$ maps $\F_T$ into itself and
$\mathcal{J}\;\mathcal{S}_T=\varphi_T$:
\[
 \mathcal{S}_T^*(x^*,x)=\varphi_T(x,x^*).
\]
Note that any $h\in \F_T$ satisfies the condition
bellow:
\begin{equation}
  \label{eq:ca.r}
\begin{array}{l}
   h(x,x^*)\geq \inner{x}{x^*}\\
   h^*(x^*,x)\geq \inner{x}{x^*}
 \end{array}\qquad
 \forall\, (x,x^*)\in X\times X^*.  
\end{equation}
What about the converse? Burachik and Svaiter proved in~\cite[Theorem
3.1]{BuSvProc03} that if a closed convex function $h$
satisfies~\eqref{eq:ca.r} in a \emph{reflexive} Banach space, then $h$
represents a maximal monotone operator and $h$ belongs to the
Fitzpatrick function of this operator.
This result has been used for ensuring maximal monotonicity in
reflexive Banach spaces~\cite{SimZaProc04,PenRelv04,PenZaProc06,BorProc06,BoGraWanSIAM06,BarBauBorNA07,BauBorWanSIAM07,BauWanSet07,BorProc07,SimJAC07}.

For the case of a non-reflexive Banach space, Marques-Alves and
Svaiter proved~\cite{MASvJCA08} that if $h$ is a convex lower
semicontinuous function in $X\times X^*$ and
\begin{equation}
  \label{eq:ca}
\begin{array}{l}
   h(x,x^*)\geq \inner{x,}{x^*},\qquad
         \forall\, (x,x^*)\in X\times X^*  \\
         h^*(x^*,x^{**})\geq \inner{x^*}{x^{**}},\qquad
         \forall\, (x^*,x^{**})\in X^*\times X^{**}
 \end{array}\qquad
\end{equation}
then again $h$ and $\mathcal{J}h$ represent a maximal monotone
operator and belong to Fitzpatrick family of this operator.
Moreover, the operator $T$ satisfies a restricted version of the Br\o
ndsted-Rockafellar property.
In particular, Marques-Alves and Svaiter proved that if $T$ is maximal monotone
and \emph{one} Fitzpatrick function of $T$ satisfies~\eqref{eq:ca}, then
$T$ satisfies the restricted Br\o ndsted-Rockafellar property.
The case of $h$ convex (but not lower semicontinuous) and
satisfying~\eqref{eq:ca} was also examined in~\cite{MASvJCA08}.

Mart\'inez-Legaz and Svaiter~\cite{LegSvSet05} 
defined (with a different notation), for $h:X\times X^*\to\BR$ and
$(x_0,x_0^*)\in X\times X^*$
\begin{equation}
  \label{eq:def.hx}
  \begin{array}{l}
    h_{(x_0,x_0^*)}:X\times X^*\to\BR,\\[.4em]
    h_{(x_0,x_0^*)}(x,x^*):=h(x+x_0,x^*+x_0^*)-[\inner{x}{x_0^*}+\inner{x_0}{x^*}
    +\inner{x_0}{x_0^*}].
  \end{array}
\end{equation}
The operation $h\mapsto h_{(x_0,x_0^*)}$ preserves many properties of
$h$, as convexity, lower semicontinuity and can be seen as the action
of the group $(X\times X^*,+)$ on $\BR^{X\times X^*}$, because
\[ \left(h_{(x_0,x_0^*)}\right)_{(x_1,x_1^*)}=h_{(x_0+x_1,x_0^*+x_1^*)}.\]
Moreover
\[
  \left(h_{(x_0,x_0^*)}\right)^*=\left(h^*\right)_{(x_0^*,x_0)},
\]
where the rightmost $x_0$ is identified with its image under the
canonical injection of $X$ into $X^{**}$.
Therefore, 
\begin{enumerate}
\item $h\geq \pi\iff h_{(x_0,x_0)}\geq \pi$,
\item $
  \left(h_{(x_0,x_0^*)}\right)^*\geq \pi_*\iff
  \left(h^*\right)_{(x_0^*,x_0)}
  \geq \pi_*$, 
\end{enumerate}
and finally,
\[ h\in\F_T\iff h_{(x_0,x_0^*)}\in \F_{T-\{(x_0,x_0^*)\}}.
\]
Marques-Alves and Svaiter work~\cite{MASvJCA08} was
heavily based on these nice properties of the map $h\mapsto
h_{(x_0,x_0^*)}$.
These authors also used the fact that
if $h$ satisfies
condition~\eqref{eq:ca}, then it also satisfies the following
\emph{auxiliary condition}: 
\begin{equation}
  \label{eq:cb}
  \inf_{(x,x^*)\in X\times X^*} h_{(x_0,x_0^*)}(x,x^*)+\frac{1}{2}\normq{x}+
\frac{1}{2}
  \normq{x^*}=0,\qquad \forall (x_0,x_0^*)\in X\times X^*. 
\end{equation}

A possible generalization of~\cite{MASvJCA08} would be to require
only the auxiliary condition \eqref{eq:cb} for one Fitzpatrick
function of a maximal monotone operator $T$ and then conclude that
this operator satisfies the restricted Br\o ndsted-Rockafellar
property.
Unfortunately, condition~\eqref{eq:cb} is not more general than
condition~\eqref{eq:ca}, as we will prove.

The class of operators studied in~\cite{MASvJCA08} is the class of
maximal {\bf M}onotone operators for which there exists a function in
Fitzpatrick family with a conjugate {\bf A}bove the duality product.
So, for the time being, we will call these operators type MA. We
will also prove that  MA condition is equivalent to NI condition. 
A maximal monotone $T:X\tos X^*$ is type (NI)~\cite{SimRanJMA96} if
\[ \inf_{(y,y^*)\in T}\inner{x^{**}-y}{x^*-y^*}\leq 0,
\forall (x^*,x^{**})\in X^*\times X^{**}.
\]
For proving this equivalence we will show that if 
 \emph{some} $h\in
\F_T$ satisfies 
condition~\eqref{eq:ca}, then
\emph{all} function in Fitzpatrick family of $T$ satisfies
condition~\eqref{eq:ca}. Observe again that, for a function in
Fitzpatrick family, $h\geq \pi$ holds by definition.

The main results of this work are the two
theorems bellow:
\begin{theorem}
  \label{th:1}
  Let $X$ be a real Banach space and $h$ be a convex function on
  $X\times X^*$. Then $h$ satisfies the condition~\cite[eq. (4)]{MASvJCA08}
    \begin{equation}
    \label{eq:t1.a}
   h\geq \pi\qquad h^*\geq \pi_*
    \end{equation}
    if, and only if, $h$ satisfies the auxiliary condition
\cite[eq. immediately bellow eq.\ (29)]{MASvJCA08},
    \begin{equation}
     \label{eq:th1.b}
      \inf_{(x,x^*)\in X\times X^*} h_{(x_0,x_0^*)}(x,x^*)+\frac{1}{2}\normq{x}+
      \frac{1}{2}
      \normq{x^*}=0,\qquad \forall (x_0,x_0^*)\in X\times X^*. 
    \end{equation}
\end{theorem}
\begin{theorem}
  \label{th:2}
  Let $X$ be a real Banach space and $T:X\tos X^*$. The following
  conditions are equivalent:
  \begin{enumerate}
  \item $T$ is type \emph{MA}, that is, $T$ is maximal monotone and there
    exists some $h\in\F_T$ such that
    $h^*\geq\pi_*$ (and  $h\geq\pi$),
  \item $T$ is maximal monotone and  \emph{all} $h\in \F_T$,
    satisfies the
    condition  $h^*\geq\pi_*$ (and  $h\geq\pi$),
  \item   $T$ is maximal monotone and \emph{some} $h\in \F_T$ satisfies the
    condition
    \[
    \inf_{(x,x^*)\in X\times X^*} h_{(x_0,x_0^*)}(x,x^*)+\frac{1}{2}\normq{x}+
    \frac{1}{2}
    \normq{x^*}=0,\qquad \forall (x_0,x_0^*)\in X\times X^*.
    \]
  \item  $T$ is maximal monotone and \emph{all} $h\in \F_T$ satisfies the
    condition
    \[
    \inf_{(x,x^*)\in X\times X^*} h_{(x_0,x_0^*)}(x,x^*)+\frac{1}{2}\normq{x}+
    \frac{1}{2}
    \normq{x^*}=0,\qquad \forall (x_0,x_0^*)\in X\times X^*.
    \]
  \item $T$ is type \emph{NI}:
    \[
    \inf_{(y,y^*)\in T}\inner{x^{**}-y}{x^*-y^*}\leq 0,\qquad \forall
    (x^*,x^{**})\in X^*\times X^{**}.
    \]
  \end{enumerate}
where $\pi$ and $\pi_*$ are the duality products in $X\times X^*$ and 
$X^*\times X^{**}$, as described in~\eqref{eq:df.pps}.
\end{theorem}

\section{Proof of the main results}
\label{sec:proof}

\begin{proof}[Proof of Theorem~\ref{th:1}]
  Let $\bar h:=\clconv h$. As $h$ is convex,
  \[ \bar h(x,x^*)=\lim\inf_{(y,y^*)\to(x,x^*)} h(y,y^*),\]
  and, for any $(x_0,x_0^*)\in X\times X^*$,
  \[
  \bar h_{(x_0,x_0^*)}(x,x^*)=\lim\inf_{(y,y^*)\to(x,x^*)}
    h_{(x_0,x_0^*)}(y,y^*).
  \] 
  As the duality product is continuous and 
  $(\clconv h)^*=h^*$, condition \eqref{eq:t1.a}
  holds for $h$ if, and only if, it holds for $\bar h$. As the norms are
  continuous (this is indeed trivial), condition \eqref{eq:th1.b}
  holds for $h$ if, and only if, it holds for $\bar h$. So, it suffices
  to prove the theorem for the case where $h$ is lower semicontinuous,
  and we assume it from now on in this proof.

  For the sake of completeness, we discuss the implication
  \eqref{eq:t1.a}$\Rightarrow$\eqref{eq:th1.b}. 
  Take $(x_0,x_0^*)\in X\times X^*$.  If condition \eqref{eq:t1.a}
  holds for $h$, then it holds for $h_{(x_0,x_0^*)}$ and using
  \cite[Theorem 3.1, eq. (12)]{MASvJCA08} we conclude that
  condition~\eqref{eq:th1.b} holds.

   For proving the implication
   \eqref{eq:th1.b}$\Rightarrow$\eqref{eq:t1.a}, first note that, for
   any $(z,z^*)\in X\times X^*$,
   \[ h_{(z,z^*)}(0,0)\geq\inf_{(x,x^*)}
   h_{(z,z^*)}(x,x^*)+\frac{1}{2}\normq{x}+ \frac{1}{2}\normq{x^*}.
   \]
   Therefore, using also \eqref{eq:th1.b} we obtain
   \[
   h(z,z^*)-\inner{z}{z^*}= h_{(z,z^*)}(0,0)\geq 0.
   \]
   Since $(z,z^*)$ is an arbitrary element of $X\times X^*$ we
   conclude that $h\geq \pi$.
  
   For proving that, under assumption  \eqref{eq:th1.b},
   $h^*\geq\pi_*$,
   take some $(y^*,y^{**})\in X^*\times X^{**}$. 
   First, use Fenchel-Young inequality to conclude that
   for any $(x,x^*), (z,z^*)\in X\times X^*$,
   \begin{align*}
     h_{(z,z^*)}(x,x^*)
     \geq&
     \inner{x}{y^*-z^*}+\inner{x^*}{y^{**}-z}-\left(h_{(z,z^*)}\right)^*(y^*-z^*,y^{**}-z).
   \end{align*}
   As $\left(h_{(z,z^*)}\right)^*=(h^*)_{(z^*,z)}$,
   \begin{align*}
     \left(h_{(z,z^*)}\right)^*(y^*-z^*,y^{**}-z)&=
     h^*(y^*,y^{**})-\inner{z}{y^*-z^*}-\inner{z^*}{y^{**}-z}-\inner{z}{z^*}\\
     &=h^*(y^*,y^{**})-\inner{y^*}{y^{**}}+\inner{y^*-z^*}{y^{**}-z}.
   \end{align*}
   Combining the two above equations we obtain
   \begin{align*}
      h_{(z,z^*)}(x,x^*)
     \geq&
     \inner{x}{y^*-z^*}+\inner{x^*}{y^{**}-z}\\
      &-\inner{y^*-z^*}{y^{**}-z}+\inner{y^*}{y^{**}}-h^*(y^*,y^{**}).
   \end{align*}
   Adding $(1/2)\normq{x}+(1/2)\normq{x^*}$ in both sides of the
   above inequality we have 
  \begin{align*}
      h_{(z,z^*)}(x,x^*)+\frac{1}{2}\normq{x}+\frac{1}{2}\normq{x^*}
     \geq&
     \inner{x}{y^*-z^*}+\inner{x^*}{y^{**}-z}
     +\frac{1}{2}\normq{x}+\frac{1}{2}\normq{x^*}
     \\
      &-\inner{y^*-z^*}{y^{**}-z}+\inner{y^*}{y^{**}}-h^*(y^*,y^{**}).
   \end{align*}
   Note that
   \[
   \inner{x}{y^*-z^*}+\frac{1}{2}\normq{x}\geq
   -\frac{1}{2}\normq{y^*-z^*},\qquad
   \inner{x^*}{y^{**}-z}+\frac{1}{2}\normq{x^*}\geq
   -\frac{1}{2}   \normq{y^{**}-z}.
   \]
   Therefore, for any $(x,x^*), (z,z^*)\in X\times X^*$,
   \begin{align*}
     h_{(z,z^*)}(x,x^*)+\frac{1}{2}\normq{x}+\frac{1}{2}\normq{x^*}
     \geq&
    -\frac{1}{2}\normq{y^*-z^*} -\frac{1}{2}   \normq{y^{**}-z}
     \\
     &-\inner{y^*-z^*}{y^{**}-z}+\inner{y^*}{y^{**}}-h^*(y^*,y^{**}).
   \end{align*}
   Using now assumption \eqref{eq:th1.b} we conclude that the infimum,
   for $(x,x^*)\in X\times X^*$, at the left hand side of the above
   inequality is $0$. Therefore, taking the infimum on $(x,x^*)\in
   X\times X^*$ at the left hand side of the above inequality and
   rearranging the resulting inequality we have
   \begin{align*}
     h^*(y^*,y^{**})-\inner{y^*}{y^{**}}\geq 
     -\frac{1}{2}\normq{y^*-z^*} -\frac{1}{2}   \normq{y^{**}-z}
     -\inner{y^*-z^*}{y^{**}-z}.
   \end{align*}
   Note that
   \[ \sup_{z^*\in X^*}  -\inner{y^*-z^*}{y^{**}-z}
   -\frac{1}{2}\normq{y^*-z^*}
   =\frac{1}{2}\normq{y^{**}-z}.
  \]
  Hence, taking the sup in $z^*\in X^*$ at the  right hand
  side of the previous inequality we obtain
  \[  h^*(y^*,y^{**})-\inner{y^*}{y^{**}}\geq 0\]
and condition~\eqref{eq:t1.a} holds.
\end{proof}

\begin{proof}[Proof of Theorem~\ref{th:2}]
  First use Theorem~\ref{th:1} to conclude that item 1 and 3 are
  equivalent. The same theorem also shows that items 2 and 4
  are equivalent.

  Now assume that item 3 holds, that is,
  for some $h\in \F_T$,
  \[
  \inf_{(x,x^*)\in X\times X^*} h_{(x_0,x_0^*)}(x,x^*)+\frac{1}{2}\normq{x}+
  \frac{1}{2}
  \normq{x^*}=0,\qquad \forall (x_0,x_0^*)\in X\times X^*.
  \]
  Take $g\in \F_T$, and $(x_0,x_0^*)\in X\times X^*$.
  First observe that, for any $(x,x^*)\in X\times X^*$, $
  g_{(x_0,x_0^*)}(x,x^*)\geq\inner{x}{x^*}$ and
  \[  g_{(x_0,x_0^*)}(x,x^*)+\frac{1}{2}\normq{x}+
  \frac{1}{2}\normq{x^*}\geq
   \inner{x}{x^*}+\frac{1}{2}\normq{x}+
   \frac{1}{2}\normq{x^*}\geq 0.
   \]
   Therefore,
   \begin{equation}
     \label{eq:aux0}
      \inf_{(x,x^*)\in X\times X^*} g_{(x_0,x_0^*)}(x,x^*)+\frac{1}{2}\normq{x}+
  \frac{1}{2}
  \normq{x^*}\geq 0. 
   \end{equation}
  As the square of the norm is coercive, there exist $M>0$ such that
  \[
  \left\{ (x,x^*)\in X\times X^*\;|\;  h_{(x_0,x_0^*)}(x,x^*)+\frac{1}{2}\normq{x}
    +  \frac{1}{2}\normq{x^*}<1
    \right\}\subset
   B_{X\times X^*}(0,M),
  \]
  where
  \[
  B_{X\times X^*}(0,M)=\left\{ (x,x^*)\in X\times X^*\;|\;
    \sqrt{\normq{x}+\normq{x^*}}<M\right\}.
  \]
  For any $\varepsilon>0$, there exists $(\tilde x,\tilde x^*)$ such
  that
  \[ \min\left\{1,\varepsilon^2\right\}>
  h_{(x_0,x_0^*)}(\tilde x,\tilde x^*)+\frac{1}{2}\normq{\tilde x}
  +  \frac{1}{2}\normq{\tilde x^*}.
  \]
  Therefore
  \begin{equation}
    \label{eq:aux1}
    \begin{array}{l}
      \varepsilon^2 >
      h_{(x_0,x_0^*)}(\tilde x,\tilde x^*)+\frac{1}{2}\normq{\tilde x}
      +\frac{1}{2}\normq{\tilde x^*}\geq
      h_{(x_0,x_0^*)}(\tilde x,\tilde x^*)
       -\inner{\tilde x}{\tilde x^*}\geq 0,\\[.5em]
     M^2\geq \normq{\tilde x}+\normq{\tilde x^*}.
    \end{array}
  \end{equation}
  In particular,
  \[  \varepsilon^2 > h_{(x_0,x_0^*)}(\tilde x,\tilde x^*)
  -\inner{\tilde x}{\tilde x^*}.
  \]
  Using now the fact that operators type MA satisfies the restricted
  Br\o ndsted-Rockafellar property~\cite[Theorem 3.4]{MASvJCA08} we conclude that there exists
  $(\bar x,\bar x^*)$ such that
  \begin{equation}
    \label{eq:aux2}
   h_{(x_0,x_0^*)}(\bar x,\bar x^*)=\inner{\bar x}{\bar x^*},\quad
  \norm{\tilde x-\bar x}<\varepsilon, \quad
  \norm{\tilde x^*-\bar x^*}<\varepsilon.   
  \end{equation}
    Therefore, 
   \[ h(\bar x+x_0,\bar x^*+x_0^*)-\inner{\bar x+x_0}{\bar x^*+x_0^*}=
   h_{(x_0,x_0^*)}(\bar x,\bar x^*)-\inner{\bar x}{\bar x^*}=0,
   \]
  and $(\bar x+x_0,\bar x^*+x_0^*)\in T$. As $g\in \F_T$,
  \[
  g(\bar x+x_0,\bar x^*+x_0^*)=\inner{\bar x+x_0}{\bar x^*+x_0^*},
  \]
  and 
  \begin{equation}
    \label{eq:aux3}
    g_{(x_0,x_0^*)}(\bar x,\bar x^*)=\inner{\bar x}{\bar x^*}.
  \end{equation}
  Using the first line of \eqref{eq:aux1} we have
  \[       \varepsilon^2 >
      h_{(x_0,x_0^*)}(\tilde x,\tilde x^*)+
      \bigg[\frac{1}{2}\normq{\tilde x}
      +\frac{1}{2}\normq{\tilde x^*}+
          \inner{\tilde x}{\tilde x^*}
      \bigg]
       -\inner{\tilde x}{\tilde x^*}\geq
\frac{1}{2}\normq{\tilde x}
      +\frac{1}{2}\normq{\tilde x^*}+
          \inner{\tilde x}{\tilde x^*}.
       \]
  Therefore,
  \begin{equation}
    \label{eq:aux4}
     \varepsilon^2>  \frac{1}{2}\normq{\tilde x}
      +\frac{1}{2}\normq{\tilde x^*}+
      \inner{\tilde x}{\tilde x^*}.
  \end{equation}
  Direct use of \eqref{eq:aux2} gives
  \begin{align*}
    \inner{\bar x}{\bar x^*}&=\inner{\tilde x}{\tilde x^*}
    +\inner{\bar x-\tilde x}{\tilde x^*}
    +\inner{\tilde x}{\bar x^*-\tilde x^*}
    +\inner{\bar x-\tilde x}{\bar x^*-\tilde x^*}\\
    &\leq \inner{\tilde x}{\tilde x^*}
    +\norm{\bar x-\tilde x}\,\norm{\tilde x^*}
    +\norm{\tilde x}\,\norm{\bar x^*-\tilde x^*}
    +\norm{\bar x-\tilde x}\,\norm{\bar x^*-\tilde x^*}\\
  &\leq \inner{\tilde x}{\tilde x^*}
    +\varepsilon [ \norm{\tilde x^*}+\norm{\tilde x}]
    +\varepsilon^2
  \end{align*}
  and
  \begin{align*}
    \normq{\bar x}+\normq{\bar x^*}&\leq
    \left( \norm{\tilde x}+\norm{\bar x-\tilde x}\right)^2
    +  \left( \norm{\tilde x^*}+\norm{\bar x^*-\tilde x^*}\right)^2\\
    &\leq \normq{\tilde x}+ \normq{\tilde x^*}
    +2\varepsilon[\norm{\tilde x}+\norm{\tilde x^*}]+2\varepsilon^2
  \end{align*}
  Combining the two above equations with~\eqref{eq:aux3} we obtain

  \[
  g_{(x_0,x_0^*)}(\bar x,\bar x^*)
  +\frac{1}{2}\normq{\bar x}
  +\frac{1}{2}\normq{\bar x^*}\leq
  \inner{\tilde x}{\tilde x^*}
  +\frac{1}{2}\normq{\tilde  x}
  +\frac{1}{2}\normq{\tilde  x^*}+2\varepsilon[\norm{\tilde x}+\norm{\tilde x^*}]+2\varepsilon^2
  \]
  Using now~\eqref{eq:aux4} and the second line of \eqref{eq:aux1} we
  conclude that
 \[
  g_{(x_0,x_0^*)}(\bar x,\bar x^*)
  +\frac{1}{2}\normq{\bar x}
  +\frac{1}{2}\normq{\bar x^*}\leq
  2\varepsilon\;M\sqrt{2}+3\varepsilon^2.
  \]
  As $\varepsilon$ is an arbitrary strictly positive number, using
  also~\eqref{eq:aux0} we conclude that
  \[
      \inf_{(x,x^*)\in X\times X^*} g_{(x_0,x_0^*)}(x,x^*)+\frac{1}{2}\normq{x}+
  \frac{1}{2}
  \normq{x^*}=0.
  \]
  Altogether, we conclude that if item 
  3 holds then item 4
  holds. 
  The converse item 4$\Rightarrow$ item 3 is trivial to verify. Hence
  item 3 and item 4 are equivalent. As item 1 is equivalent to 3 and
  item 2 is equivalent to 4, we conclude that items 1,2,3  and 4 are
  equivalent.

Now we will deal with item 5. First suppose that item 2 holds. Since 
$\mathcal{S}_T\in \F_T$ 
\[
 (\mathcal{S}_T)^*\geq \pi_*.
\]
As has already been observed, for any proper function $h$ it holds that $(\mbox{cl\,conv}\,h)^*=h^*$.
Therefore 
\[
 (\mathcal{S}_T)^*=(\pi+\delta_T)^*\geq \pi_*,
\]
that is 
\begin{equation}\label{eq20}
\sup_{(y,y^*)\in T}
 \inner{y}{x^*}+\inner{y^*}{x^{**}}-\inner{y}{y^*}\geq
 \inner{x^*}{x^{**}}, \forall
    (x^*,x^{**})\in X^*\times X^{**}
\end{equation}
After some algebraic manipulations we conclude
that~\eqref{eq20} is equivalent to
\[
    \inf_{(y,y^*)\in T}\inner{x^{**}-y}{x^*-y^*}\leq 0,\qquad \forall
    (x^*,x^{**})\in X^*\times X^{**},
\]
that is, $T$ is type (NI).
If item 5 holds, by the same reasoning we conclude that~\eqref{eq20} holds and therefore
$(\mathcal{S}_T)^*\geq \pi_*$. As $\mathcal{S}_T\in \F_T$, we conclude that item 5 $\Rightarrow$ item 1. As has been proved previously item 1 $\Rightarrow$ item 2.

\end{proof}

\bibliographystyle{plain}

\end{document}